\documentclass[11pt]{article}

\usepackage{amsmath}
\usepackage{amsfonts}

\newtheorem{tw}{Theorem}[section]
\newtheorem{lem}[tw]{Lemma}
\newtheorem{cor}[tw]{Corollary}
\newtheorem{rem}[tw]{Remark}
\newtheorem{df}[tw]{Definition}
\newtheorem{nt}[tw]{Notation}
\newcommand{\bl}{\begin{lem}}
\newcommand{\el}{\end{lem}}
\newcommand{\bn}{\begin{nt}\rm}
\newcommand{\en}{\end{nt}}
\newcommand{\bd}{\begin{df}}
\newcommand{\ed}{\end{df}}
\newcommand{\br}{\begin{rem}\rm}
\newcommand{\er}{\end{rem}}
\newcommand{\bt}{\begin{tw}}
\newcommand{\et}{\end{tw}}
\newcommand{\bc}{\begin{cor}}
\newcommand{\ec}{\end{cor}}
\newcommand{\N}{{\mathbb{N}}}
\newcommand{\Z}{{\mathbb{Z}}}
\newcommand{\R}{{\mathbb{R}}}
\newcommand{\supp}{{\rm supp\:}}
\newcommand{\proof}{\noindent{\bf Proof.} }
\newcommand{\ca}{{\mathcal{A}}}
\newcommand{\cn}{{\mathcal{N}}}
\newcommand{\nf}{{\N^{<\infty}}}
\newcommand{\cm}{\mathcal{M}}
\newcommand{\kr}{{K}_{q,r}}
\newcommand{\con}{\subseteq}
\newcommand{\conj}[2]{ \{ {#1}\,:\,{#2} \} }
\newcommand{\al}{\alpha}
\newcommand{\vep}{\varepsilon}

\title{Stabilization of Tsirelson-type norms on $\ell_p$ spaces}
\author{Anna Pelczar}
\begin{document}

\maketitle
\begin{abstract}
We consider classical Tsirelson-type norms of $T[\ca_n,\theta]$
and their modified versions on $\ell_p$ spaces. We show that for
any $1<p<\infty$ there is a constant $\lambda_p$ such that
considered Tsirelson-type norms do not $\lambda_p$-distort any of
subspaces of $\ell_p$.
\end{abstract}

\section{Preliminaries}

The Tsirelson space $T$ was the base for further constructions of Banach spaces, which
solved crucial problems in the theory of Banach spaces, as Schlumprecht space - the first
space known to be arbitrary distortable, and further on the hereditarily indecomposable
Gowers-Maurey space. Tsirelson-type norms provide important examples of Banach spaces as
well as a uniform approach to the classical spaces and the new ones - by defining norms
implicitly, as solutions to certain equalities. The mixed and modified mixed
Tsirelson-type norms were studied in various context, with respect to their distortion
and asymptotic properties (cf. \cite{ao,ad,ad2,adm,adkm}).

The modified Tsirelson norm on $T$ was introduced in \cite{j} and later proved to be
equivalent to the Tsirelson norm in \cite{co} and in \cite{b} with constant 2. Other
equivalent norms $\|\cdot\|_n$ on $T$  - of spaces $T[\mathcal{S}_n,\frac{1}{2^n}]$ -
were studied in \cite{ot} in context of the question of arbitrary distortability of $T$.
It was shown that there is a universal constant $K$, such that the norms $\|\cdot\|_n$ do
not $K$-distort any infinite dimensional subspace of $T$. Their modified versions
considered in \cite{m} appeared to be 3-equivalent to the original versions.

The equivalence of certain classical Tsirelson-type norms and
their modified versions was shown in \cite{b}. In the same paper
the norms of $T[\ca_n,\theta]$ isomorphic to $\ell_p$ were proved
to be $\theta^2$-equivalent to the classical $\|\cdot\|_p$ norms
on $\ell_p$. General Tsirelson-type norms appeared in \cite{ad},
where equivalence between classical $\ell_p$ norm and
Tsirelson-type norms of $T[\ca_n,\theta]$ was shown by means of
tree analysis of norming vectors. The fact, applied in our paper,
that classical and modified Tsirelson-type norms on $\ell_p$
spaces are 3-equivalent follows from \cite{pt}. Let us recall that
in the mixed case original and modified version of Tsirelson-type
norms define non-isomorphic spaces, as in the case of Schlumprecht
space (cf. \cite{adkm}).

E.Odell and T.Schlumprecht solved in \cite{os} the famous distortion problem showing that
the spaces $\ell_p$, for $1<p<\infty$, are arbitrary distortable. They have shown in fact
that these spaces are biorthogonally distortable, transferring the so-called biorthogonal
system from the Schlumprecht space. The question of norms on $\ell_p$ arbitrary
distorting the original norm, defined only by means of $\ell_p$, remains open. The
obvious candidates to be studied in this context are the Tsirelson-type norms.

We show in this paper that for any $1<p<\infty$ there is a constant $\lambda_p$ such that
the Tsirelson-type norms of $T[\ca_n,\theta]$ do not $\lambda_p$-distort any of infinite
dimensional subspaces of the $\ell_p$ space. In fact we prove stabilization of modified
Tsirelson-type norms on $\ell_p$ (Theorem \ref{maintheorem}), which are 3-equivalent to
the original Tsirelson-type norm (Corollary \ref{modified}) by the reasoning of
\cite{pt}.

The author would like to thank Jordi Lopez-Abad for valuable
remarks and simplifying certain proofs.

\

We recall first the standard notation. By $c_{00}$ we denote the space of real sequences
which are eventually zero, endowed with the supremum norm $\|\cdot\|_{\infty}$; by
$\ell_p$ ($1<p<\infty$) the space of $p$-summable real sequences with the canonical norm
$\|\cdot\|_p$ given by the formula

$$\|x\|_p=\left(\sum_{i=1}^{\infty}|x(i)|^p\right)^{1/p} \mbox{ for any } x=(x(i))_i\in\ell_p$$

By $(e_n)$ we denote the unit vectors basis. As usual we put
$B_{\ell_p}=\conj{x\in\ell_p}{\| x\|_p\le 1}$.

For $1<p<\infty$ we have $(\ell_p)^* \cong\ell_q$, where $\frac{1}{p}+\frac{1}{q}=1$. We
put
$$\langle x,y\rangle=\sum_i x(i)y(i) \mbox{ for any } x=(x(i))_i\in\ell_p \mbox{ and } y=(y(i))_i\in\ell_q$$

For any sets $I,J\con \N$ we write $I<J$ if $\max I <\min J$ and for any vectors $x,y\in
c_{00}$ we write $x<y$ if $\supp x<\supp y$. A sequence $(x_n)\con c_{00}$ is called a
block sequence provided $x_1<x_2<\dots$. Given a block sequence $(x_n)$ by $[x_n]$ we
denote the vector space spanned by $(x_n)$.

Given any $x\in c_{00}$ and $E\con \N$ by $Ex$ or $x_E$ denote the restriction of $x$ to $E$, ie. $Ex(i)=x(i)$
if $i\in E$ and $Ex(i)=0$ otherwise.

Finally, we say that a set $K\con c_{00}$ is closed (or invariant) under

\noindent (a) restriction, if for any $x\in K$ and $E\subset \N$ also $Ex\in K$,

\noindent (b) spreading, if for any $x=\sum a_n e_n \in K$ and any strictly increasing
function $\phi:\N\to \N$ we have $\sum a_n e_{\phi(n)}\in K$,

\noindent (c) permutation, if for any $x=\sum a_n e_n \in K$ and any permutation
$\sigma:\N\to \N$ we have $\sum a_n e_{\sigma(n)}\in K$.

\section{Classical and modified Tsirelson-type norms}

We recall briefly the construction of Tsirelson-type norm, denoted here by
$\|\cdot\|_{p,r}$, and modified Tsirelson-type norm, denoted by $|\cdot|_{p,r}$ (cf.
\cite{ad}, \cite{adkm}).

\bd

Fix $1<p,q<\infty$ with $\frac{1}{p}+\frac{1}{q}=1$ and $r\in\N$.

Define norms $\|\cdot\|_{p,r}$ and $|\cdot|_{p,r}$ on $c_{00}$ as the unique norms satisfying the following
implicit equations for any $x\in c_{00}$:

$$\| x\|_{p,r}=\max\left\{\|x\|_{\infty},\frac{1}{\sqrt[q]{r}}\sup\sum_{i=1}^r\|E_ix\|_{p,r}\right\}$$

where the supremum is taken over all $r$-tuples of sets $E_1,\dots,E_r\con\N$ which satisfy $E_1<\dots<E_r$,

$$|x|_{p,r}=\max\left\{\|x\|_{\infty},\frac{1}{\sqrt[q]{r}}\sup\sum_{i=1}^r|F_ix|_{p,r}\right\}$$

where the supremum is taken over all $r$-tuples of sets $F_1,\dots,F_r\con\N$, which are pairwise disjoint.

\ed

\br

\begin{enumerate}

\item[(a)] Basic unit vectors $(e_n)$ form an 1-unconditional and 1-subsymmetric basis of
$c_{00}$ endowed with $\|\cdot\|_{p,r}$ and an 1-unconditional and 1-symmetric basis of
$c_{00}$ endowed with $|\cdot |_{p,r}$.

\item[(b)] The completion of $c_{00}$ endowed with the norm $\|\cdot\|_{p,r}$, $r>1$, is
isomorphic to $\ell_p$ \cite{ad}. As we have $\|\cdot\|_{p,r}\le |\cdot
|_{p,r}\le\|\cdot\|_p$ it follows that the completion of $c_{00}$ endowed with the norm
$|\cdot |_{p,r}$, $r>1$, is also isomorphic to $\ell_p$.

\end{enumerate}
\er

The norms introduced above can be defined alternatively by their norming sets presented
below.

\bd

Fix $1<q<\infty$ and $r\in\N$.

Let $K_{q,r}$ be the smallest set in $c_{00}$ which contains vectors $(\pm e_n)$ and
satisfies the following:
$$z_1,\dots, z_l \in \kr,\;\; l\le r,\;\; z_1<\dots<z_l
\;\;\Longrightarrow\;\; \frac{1}{\sqrt[q]{r}}(z_1+\dots+z_l)\in \kr.
$$

Let $\kr^{\cm}$ be the smallest set in $c_{00}$ which contains
vectors $(\pm e_n)$ and satisfies the following:
$$y_1,\dots, y_l\in \kr^{\cm},\;\; l\le r,\;\; \supp y_i\cap\supp y_j=\emptyset,\;\; i\neq
j\;\;\Longrightarrow\;\; \frac{1}{\sqrt[q]{r}}(y_1+\dots+y_l)\in
\kr^{\cm}.$$

\ed

\br \begin{enumerate}

\item[(a)] By definition, in particular the minimality of
considered sets, we have $\kr\con\kr^{\cm}\con B_{\ell_q}$.

\item[(b)] By definition sets $\kr$ and $\kr^{\cm}$ are closed
under restriction and spreading. The set $\kr^{\cm}$ is a
"symmetrized" version of $\kr$ closed under permutations.

\item[(c)] A standard reasoning proves that for any $r\in\N$, $1<p,q<\infty$ with
$\frac{1}{p}+\frac{1}{q}=1$ we have
 \begin{align*}
\| x\|_{p,r}=&\;\sup\conj{\langle x,z\rangle}{z\in \kr}\\
 |x |_{p,r}=&\;\sup\conj{\langle x,y \rangle}{y\in \kr^{\cm}}
 \end{align*}

\item[(d)] As for $r=1$ clearly $K_{p,1}=\{\pm e_n\}$, we will omit this case in the rest
of the paper.

\end{enumerate}
\er

\bd Given $\al >0$ and $1<q<\infty$ put
 \begin{align*}
C_{\al}= &\;\conj{\pm\al^j}{j\in \Z }\cup\{0\}\\
\cn_{\al}=&\; \conj{ x\in c_{00}}{x(i)\in C_{\al}, i\in\N}\\
\cn_{\al}^{(q)}=&\; \cn_{\al}\cap B_{\ell_q}
\end{align*}
 \ed

In the rest of the paper we shall need the following
characterization of the set $\kr^{\cm}$:

\bl\label{form} For any $r\in\N$, $1<q<\infty$, $t=\sqrt[q]{r}$,
we have $\kr^{\cm}=\cn_{t}^{(q)}$. \el

\proof Obviously $\kr^{\cm}\con\cn_{t}^{(q)}$.

Suppose now that $x\in \cn_{t}^{(q)}$. It is easy to see that
there is some $y\in \cn_{t}^{(q)}$ such that $x<y$ and $\|
x+y\|_q=1$. Since $\kr^{\cm}$ is closed under restrictions, we may
assume that $\|x\|_q=1$. Since both $\kr^{\cm}$ and
$\cn_{t}^{(q)}$ are closed under permutations, we may assume that
$|x(1)|\ge |x(2)|\ge\dots$. The proof goes by induction on
$$n(x)=\min\conj{n}{
\mbox{there is some $j\in \supp x$ with }|x(j)|=t^{-n} }.$$
 If $n(x)=0$, then the result is clear. Suppose now that $n(x)>0$.

\

\noindent \textbf{Claim} \emph{There is a block sequence $(x_i)_{i=1}^r\con
\cn_{t}^{(q)}$ such that $x=x_1+\dots+x_r$ and $\|x_i\|_q=t^{-1}$, for any $1\le i\le
r$.}

\

\noindent \textbf{Proof of Claim.} The proof goes by induction on
$$m(x)=\max\conj{n}{ \mbox{there is some $j\in \supp x$ with }|x(j)|=t^{-n} }$$

Let $I= \conj{j\in \supp x}{|x(j)|> t^{-m(x)}}$ and $J=\supp x\setminus I$. Notice that
$I$ is an initial part of $\supp x$. Denote by $x_I$ and $x_J$ the projections of $x$ on
$I$ and $J$ respectively. Then we have that
$$1=\|x\|_q^q=\|x_I\|_q^q+ \|x_J\|_q^q=\|x_I\|_q^q+\frac{|J|}{r^{m(x)}}=
\frac{l}{r^{m(x)-1}}+\frac{|J|}{r^{m(x)}}$$
 for some integer $l\in \N$. Since $m(x)\geq
n(x)\ge 1$ it follows that $|J|=k r$ for some integer $k\in \N$. Divide $J$ into $k$
disjoint pieces $(J_i)_{i=1}^k$, with $J_1<\dots<J_k$ and $|J_i|=r$ ($1\le i\le k$). Now
pick $n_i\in J_i$ ($1\le i\le k$) and set
$$ y=x_I+t^{1-m(x)}\sum_{i=1}^k e_{n_i}$$

It is clear that $m(y)=m(x)-1$, $\|y\|_q=\|x\|_q=1$, and $|y(1)|\ge |y(2)|\ge \dots$,
hence by inductive hypothesis there is a decomposition $y=y_1+\dots +y_r$ into a sum of a
block sequence with $\|y_i\|_q=t^{-1}$ for any $1\le i\le r$.

Define $F:\supp y\to \kr^{\cm}$ by $F(j)=e_j$ if $j\in I$, and
$F(n_i)=t^{-1}\sum_{n\in J_i}e_n$, for $1\le i\le k$. It is clear
that $\|F(j)\|_q=1$ for any $j\in\supp y$ and that $F(j)<F(j')$
for any $j<j'$. For any $1\le i\le r$ define
$$x_i=\sum_{j\in \supp y_i} y_i(j)F(j)$$

By the previous observations, we obtain that $\|x_i\|_q=\|y_i\|_q=t^{-1}$ for any $1\le
i\le r$ and $(x_i)$ is a block sequence, therefore we have the decomposition
$x=x_1+\dots+ x_r$. \hfill $\square$

\

Now we continue the proof of the Lemma \ref{form}. Take the
decomposition $x=x_1+\dots +x_r$ as in the Claim. Then $\|t
x_i\|_q=1$, and $n(x_i)\le n(x)-1$ for any $1\le i\le r$, hence by
inductive hypothesis we have that $(t x_i)_i\con \kr^{\cm}$. Hence
$x=t^{-1}(tx_1+\dots +tx_r)\in \kr^{\cm} $. \hfill $\square$

\section{Equivalence of $|\cdot |_{p,r}$ and $\|\cdot\|_{p,r}$ norms}

The fact that $\|\cdot\|_{p,r}$ and $|\cdot |_{p,r}$ are
3-equivalent follows immediately from results in \cite{pt}. We
recall the reasoning from this preprint for the sake of
completeness.

\

First we introduce some notation. Let $\nf$ denote the set of finite sequences of $\N$.
For any $m= (m(1),\dots,m(n))\in\nf$ and $k\in\Z$ put $m+k{\bf 1}=(m(1)+k,\dots,m(n)+k)$.
Given any $m,l\in\nf$, $m=(m(1),\dots, m(n))$, $l=(l(1),\dots, l(j))$ put $m^{\frown}
l=(m(1),\dots, m(n),l(1),\dots, l(j))$.

\

Fix $1<q<\infty$ and $r\in\N$.

Define the function $\Phi:\nf\rightarrow (0,\infty)$ in the following way:
\begin{align*}
\Phi (m(1))= & \; r^{-m(1)} \\
\Phi (m(1),m(2))= & \; r^{-m(1)}+r^{-m(2)}\\
\Phi (m(1),\dots,m(n))= & \; r^{-m(1)}+2\sum_{i=2}^{n-1}r^{-m(i)}+r^{-m(n)}, \;\; n>2
\end{align*}

Notice that the function $\Phi$ has the following property:
$$\Phi (m(1),\dots, m(n))=\Phi (m(1),\dots, m(i))+\Phi(m(i),\dots, m(n))\;\;\; {\rm for}\;\;\; 1<i<n$$

Put $t=\sqrt[q]{r}$ and define the function
$$ V:\nf\ni (m(1),\dots,m(n))\mapsto (t^{-m(1)},\dots, t^{-m(n)},0\dots)\in c_{00}$$

\bt\label{jacek}{\rm\cite{pt}} With the above notation let the
sequence $m\in\nf$ satisfy $\Phi(m)\le 1$. Then $ V(m)\in \kr$.

\et

\noindent\textbf{Proof} goes by induction on the length $n$ of the sequence.

For $n=1$ the assertion holds true since $(t^{-m(1)},0,\dots)\in \kr$ for any $m(1)\in\N$.

Fix $n\in\N$ and assume that the Theorem holds true for any sequence of integers of length less or equal to $n$
and pick some sequence $m=(m(1),\dots, m(n+1))$ of integers.

Let us first notice that we can consider only the case $r^{-1}<\Phi (m)\le 1$. Indeed,
for any $m\in\nf$ pick $k\in\N$ such that $r^{-k-1}<\Phi (m)\le r^{-k}$. Notice that
$r^{-1}<\Phi (m-k{\bf 1})\le 1$. If $ V(m-k{\bf 1})\in \kr$, then also
$$ V(m)=t^{-k} V(m-k{\bf 1})\in \kr$$

Let now $r^{-1}<\Phi (m)\le 1$. Put $k_0=1$ and define inductively $k_1<\dots<k_l=n+2$ as
$$k_{i+1}=\max\{k\in\{k_i,\dots,n+2\}:\;\Phi
(m(k_i),\dots,m(k-1))\le r^{-1}\},\;\;i\ge 1$$

Since $r^{-1}<\Phi (m)$ we have $k_1\le n+1$.

We will show that $l\le r$. Assume that $l>1$. By the definition of $k_{i+1}$ we have
$$\Phi (m(k_i),\dots, m(k_{i+1}))>r^{-1}\;\;\; {\rm for}\;\;\; 0\le i\le l-2$$

By the addition rule for $\Phi$ we have
$$1\ge\Phi (m)=\sum_{i=0}^{l-1}\Phi (m(k_i),\dots,m(k_{i+1}))>\sum_{i=0}^{l-2}\frac{1}{r}=\frac{l-1}{r}$$

which implies that $l-1<r$, hence $l\le r$.

Define sequences $m_1,\dots, m_l$ by
$$m_i=(m(k_i),\dots, m(k_{i+1}-1))\;\;\; {\rm for}\;\;\; 0\le i\le l-1$$

Since $k_l=n+2$ we have $m={m_1}^{\frown}\dots^{\frown} m_l$. By construction the length
of $m_i$ is less or equal $n$ and $\Phi (m_i-{\bf 1})\le 1$ for any $1\le i\le l$. By the
inductive hypothesis $ V(m_i-{\bf 1})\in \kr$ for any $1\le i\le l$. Notice that
$$ V(m)= V({m_1}^{\frown}\dots^{\frown} m_l)=
t^{-1}V((m_1-{\bf 1})^{\frown}\dots^{\frown}(m_l-{\bf 1}))=t^{-1}(v_1+\dots+v_l),$$
 where $(v_1,\dots, v_l)$ is a block sequence of properly shifted vectors $V(m_1-{\bf 1}),\dots,
V(m_l-{\bf 1})$. By the definition of the set $\kr$ and its invariance under spreading it
follows that $ V(m)\in \kr$.\hfill $\square$

\bc \label{modified} Fix $1<p<\infty$ and $r\in\N$. Then
$$\frac{1}{3} |\cdot |_{p,r}\le\|\cdot\|_{p,r}\le |\cdot |_{p,r}$$ \ec

\proof Take $1<q<\infty$ with $\frac{1}{p}+\frac{1}{q}=1$. Notice that
$$\kr^{\cm}\cap 2^{-1/q}B_{\ell_q}\con \kr$$ Indeed, it follows immediately
from Theorem \ref{jacek}, since sets $\kr$ and $\kr^{\cm}$ are
invariant under permutation and for any $m\in\nf$ we have $\Phi
(m)\le 2\| V(m)\|_q^q$.

Take arbitrary $y\in \kr^{\cm}$. If for some $i\in\N$ we have
$y(i)=1$, then $y\in \kr$. If for all $i\in\N$ we have $y(i)<1$,
then $y=y_1+y_2+y_3$ for some $y_1<y_2<y_3$ with $\|y_j\|_q^q\le
\frac{1}{2}$ for $j=1,2,3$. By the above $y_1,y_2,y_3\in \kr$.

Now fix $x\in c_{00}$ and compute
 $$|x |_{p,r}=\sup_{y\in \kr^{\cm}} |\langle x,y\rangle|\le
\sup_{y_1,y_2,y_3\in \kr}|\langle x,y_1+y_2+y_3\rangle| \le
 3\sup_{z\in \kr}|\langle x,z\rangle|= 3\|x\|_{p,r}$$
 which proves the first inequality. The second inequality is obvious. \hfill $\square$

\section{Stabilization of $|\cdot |_{p,r}$ norms on $\ell_p$}

Now we present the main theorem of this paper

\bt\label{maintheorem} Fix $1<p<\infty$ and $r\in\N$, $r>1$. Every infinite dimensional
subspace $X\con \ell_p$ has an infinite dimensional subspace $Y\con X$ such that for any
$x\in Y$
$$ 4^{-6}\|x\|_p\le (\log_2 r)^{1/p}\;| x |_{p,r}\le 3\cdot 4^7 (p+q) \|x\|_p$$
 where $1<q<\infty$ satisfies $\frac{1}{p}+\frac{1}{q}=1$. \et

Thus the stabilization constant $\lambda_p$ of the norm
$|\cdot|_{p,r}$ on $\ell_p$ is not greater than $3\cdot
4^{13}(p+q)$.

\

Throughout this section we will use the following

\

\noindent \textbf{Notation}. Fix $1<p,q<\infty$ with $\frac{1}{p}+\frac{1}{q}=1$ and
$r\in\N$, $r>1$. Put $t=\sqrt[q]{r}$ and $s=\sqrt[p]{r}$.

\

Before proceeding to the proof of the main theorem we give two straightforward lemmas
comparing $\|\cdot\|_p$ and $|\cdot |_{p,r}$ on vectors from $\cn_s$.

\bl \label{comparing1} For any $x\in \cn_s^{(p)}$ we have $\|x\|_p^p\le |x |_{p,r}$\el

\proof Given $x=(x(i))\in \cn_s^{(p)}$ consider a vector $y=(y(i))\in\cn_t^{(q)}$ defined
by
$$y(i)={\rm sign}\; (x(i))|x(i)|^{p/q},\;\; i\in\N.$$
Observe that $\|x\|_p^p= \langle x,y \rangle$, which by Lemma \ref{form} ends the proof.
\hfill $\square$

\bl\label{comparing2} For any $x\in\cn_s$ we have $|x|_{p,r}\le 1+\frac{\|x\|_p^p}{t}$. \el

\proof Take $x\in\cn_s$ and any vector $y\in \cn_t^{(q)}$. Put $I=\conj{i\in\N}{
|x(i)|^p\le |y(i)|^q}$, $J=\conj{i\in\N}{|x(i)|^p>|y(i)|^q}=\conj{i\in\N}{|x(i)|^p\ge
r|y(i)|^q}$. Using the relation $\frac{1}{p}+\frac{1}{q}=1$ we obtain
$$|\langle x,y \rangle|=\sum |x(i)y(i)|
\le\sum_{i\in I}|y(i)|^q+\sum_{i\in J}\frac{|x(i)|^p}{t}\le 1+\frac{\|x\|_p^p}{t}$$
\hfill $\square$

We recall that for any $a\in\R$ by $[a]$ we denote the biggest integer smaller than $a$.

\

\noindent\textbf{Notation}. Let $M=[\log_2 r]$ and take $\al>0$ such that $\al^M=s$.

\

By definition we have $\al^{\frac{1}{2}\log_2 r}\le \sqrt[p]{r}\le \al^{\log_2 r}$.
Applying $\log_2$ to those inequalities we obtain the following

\br\label{al} $\sqrt[p]{2}\le \al \le \sqrt[p]{4}$. \er

\noindent\textbf{Notation}. For any vector $x\in\cn_{\al}$ and any $m\in\N$ we put
$$J_{m,x}=\conj{i\in\supp
x}{\al^{-m}x(i)\in C_s},\;\;\;\; J_mx=J_{m,x}x.$$ Notice that $J_{kM+m}x=J_mx$ for any
$k,m\in\N$. Using this notation we can write any vector $x\in\cn_{\al}$ as a sum
$x=J_0x+\dots+J_{M-1}x$ of vectors with disjoint supports such that
$\al^{-m}J_mx\in\cn_s$ for any $0\le m\le M-1$.

\

The proof of the main theorem uses two lemmas - Lemma \ref{stab}
showing stabilization of $|\cdot|_{p,r}$ norms on subspaces of a
certain form, and Lemma \ref{approxim} implying the saturation of
$\ell_p$ by subspaces close to those used in Lemma \ref{stab}. In
both Lemmas we will work only on vectors from $\cn_{\al}$. It is
sufficient, as for any seminormalized block sequence $(v_n)$ in
$\ell_p$ there is an $\al$-equivalent block sequence
$(x_n)\con\cn_\al$.

\bl\label{stab} Let $(x_n)_n\con\cn_{\al}$ be a block sequence such that $\al^{-3}\le\|J_mx_n\|_p\le 1$, for
$n,m\in\N$. Then for any $a_1,\dots,a_N\in\R$ we have
$$\al^{-4p}\|a_1x_1+\dots+a_Nx_N\|_p\le \sqrt[p]{M}|
a_1x_1+\dots+a_Nx_N |_{p,r}\le 6(p+q)\al^4 \|a_1x_1+\dots+a_Nx_N\|_p.$$ \el

\noindent \textbf{Proof} is based on the fact that - roughly speaking - for any block $x$
of vectors with parts $J_0,\dots,J_{M-1}$ of almost equal $\|\cdot\|_p$ norm the norm
$|x|_{p,r}$ is determined by the $|J_0x|_{p,r}$ norm.

\

Take $(x_n)$ as in Lemma. Notice that in particular $\al^{-3}\sqrt[p]{M}\le\| x_n\|_p\le
\sqrt[p]{M}$ for any $n\in\N$. We can consider only scalars $a_1,\dots, a_N\in\R$ with
$a_1,\dots,a_N>0$ and $a_1^p+\dots+a_N^p=1$. Then
$$\al^{-3}\sqrt[p]{M}\le\|\sum a_nx_n\|_p\le \sqrt[p]{M}$$

We approximate $\sum a_nx_n$ by vectors from $\cn_{\al}$. For any $1\le n\le N$ pick
$k_n\in\N$ such that $ \al^{-k_n}\le a_n<{\al^{-k_n+1}}$. Put $x=\sum_n{\al^{-k_n}}x_n\in
\cn_\al$. Notice that, for $\|\cdot\|$ denoting either $\|\cdot\|_p$ or $|\cdot|_{p,r}$,
we have
$$ \al^{-1}\|\sum a_n x_n\|\le \|x\|\le
\|\sum a_n x_n\|$$
 Observe that for any $m\in\N$ we have
$$J_mx=\sum_{n=1}^N \al^{-k_n}J_{m+k_n}x_n$$
 As $\|J_mx\|_p^p=\sum
{\al^{-pk_n}}\|J_{m+k_n}x_n\|_p^p$ by the assumptions on $(x_n)$ and the choice of
$(k_n)$ we have
$$\al^{-4}\le \|J_mx\|_p\le 1 ,\;\; m\in\N$$

\

\noindent \textsc{Left estimate.} By Lemma \ref{comparing1}, since $J_0x\in \cn_s^{(p)}$, we have $\al^{-4p}\le
|J_0x |_{p,r}$. Therefore we have
$$|\sum a_nx_n |_{p,r}\ge
|x |_{p,r}\ge |J_0x |_{p,r} \ge\al^{-4p}\ge\al^{-4p}\frac{1}{\sqrt[p]{M}}\|\sum
a_nx_n\|_p$$

\

\noindent \textsc{Right estimate}. Since $\al^{M-m}J_m x\in \cn_{s}$, by Lemma \ref{comparing2} we obtain that

$$|J_m x|_{p,r} \le \al^{m-M}+
\frac{\al^{(M-m)(p-1)}}{t}\|J_mx\|_p^p \le \al^{m-M}+ \frac{ 1}{\al^{m(p-1)}},$$

Therefore, using Remark \ref{al} and the fact that $2^{a}-1\ge a\ln 2$ for any $a>0$, we
have
\begin{align*}
\al^{-1}|\sum a_n x_n|_{p,r}\le |x|_{p,r}\le &
\sum_{m=1}^M\frac1{\al^m}+\sum_{m=0}^{M-1}\frac{ 1}{\al^{m(p-1)}} =
\frac{s-1}{s}\frac{1}{\al-1}+ \frac{s^{p-1}-1}{s^{p-1}}\frac{\al^{p-1}}{\al^{p-1}-1}\le \\
\le\frac{1}{\sqrt[p]{2}-1}+\frac{\al^{p-1}}{\sqrt[q]{2}-1}\le &
\frac{p}{\ln 2}+ \frac{\al^{p-1} q}{\ln 2}\le 6(p+q)\le
6(p+q)\frac{\al^3}{ \sqrt[p]{M}}\|\sum a_nx_n\|_p
\end{align*}

Hence
$$|\sum a_n x_n|_{p,r}\le 6(p+q)\frac{\al^4}{
\sqrt[p]{M}}\|\sum a_nx_n\|_p$$ \hfill $\square$

\bl\label{approxim} Let $(x_n)\con \cn_\al$ be a block sequence with $\al^{-2}\le
\|x_n\|_p\le \al^{-1}$ for any $n\in \N$. Then there is a block sequence $(y_n)\con
\cn_{\al}$ of $(x_n)$ such that $\al^{-3}\le \|J_m y_n\|_p\le 1$ for any $n,m\in\N$.

\el

\proof In order to prove the Lemma it will be sufficient to find one vector $y$ with the
property described above. Without loss of generality, passing to a subsequence if
necessary, we may assume that for a fixed $\vep>0$ there are scalars $b_0,\dots,b_{M-1}$
such that
$$0\le b_m-\|J_mx_n\|_p^p<\vep,\;\;\; 0\le m< M,\;\; n\in \N$$

Observe that $\al^{-2p}\le \sum_{m=0}^{M-1}b_m < \al^{-p}+ M\vep $. For technical
reasons, we define also $b_{M+m}=b_m$ for any $0\le m<M$. Take $a\in\R$ large enough so
that $[b]/b\ge 1-\vep$ for any $b\ge a$. Now fix $l\in\N$ such that $\al^{Ml}\ge a$. Now
for $0\le k,m<M$ consider the following averages:
\begin{align*}
y_k= & \frac{1}{\al^{Ml+k}}(x_{r^{l+1}k}+\dots + x_{r^{l+1}k+[\al^{(Ml+k)p}]-1})\\
 c_k^m= & \frac{[\al^{(Ml+k)p}]}{\al^{(Ml+k)p}} b_{m+k}.
\end{align*}

It is straightforward to check that $(y_k)$ is a block sequence. The sequences defined
above have the following properties:

\begin{enumerate}

\item[(a)] $J_m y_k=\al^{-(Ml+k)}J_{m+k}(x_{r^{l+1}k}+\dots +
x_{r^{l+1}k+[\al^{(Ml+k)p}]-1})$ for any $0\le k,m<M$,

\item[(b)] $0\le c_k^m-\|J_my_k\|_p^p<\vep$ for any $0\le k,m<M$,

\item[(c)] $(1-\vep) \al^{-2p} \le \sum_{k=0}^{M-1} c_k^m< \al^{-p}+ M\vep$ for any $0\le
m<M$.

Set $y=y_0+\dots +y_{M-1}$. Then

\item[(d)] $ (1-\vep)\al^{-2p} -M\vep < \|J_m y\|_p^p< \al^{-p}+M\vep$ for any $0\le
m<M$.

\end{enumerate}

Choosing sufficiently small $\vep$ we obtain the desired result. \hfill $\square$

\

\noindent \textbf{Proof of Theorem \ref{maintheorem}.} Take any infinite dimensional
subspace $X$ of $\ell_p$.

Pick any sequence $(v_n)_n\con X$ converging weakly to zero with $\| v_n\|_p=\al^{-3/2}$,
$n\in\N$. By a well-known procedure applied simultaneously to norms $\|\cdot\|_p$ and
$|\cdot |_{p,r}$ we pick for any $\delta>0$ a block sequence $(u_n)$ with
$\|u_n\|_p=\al^{-3/2}$, $n\in\N$, which is $(1+\delta)$-equivalent to some subsequence
$(v_{l_n})_n$ in both $\|\cdot\|_p$ and $|\cdot |_{p,r}$ norms.

Approximate vectors $(u_n)_n$ by vectors from $\cn_\al$: for any $n\in\N$ and $i\in\supp
u_n$ pick $k_n(i)\in\N$ such that
$$\al^{-1/2}|u_n(i)|\le \al^{-k_n(i)}\le\al^{1/2}|u_n(i)|$$
and define $(x_n)$ by conditions: $\supp x_n=\supp u_n$, $|x_n(i)|=\al^{-k_n(i)}$ and
sign $x_n(i)$= sign $u_n(i)$ for any $i\in\supp x_n$ and $n\in\N$.

The sequence $(x_n)\con\cn_{\al}$ is a block sequence with
$\al^{-2}\le\|x_n\|_p\le\al^{-1}$ for any $n\in\N$. By Lemma \ref{approxim} there is a
block sequence $(y_n)$ of $(x_n)$ satisfying assumptions of Lemma \ref{stab}.

Notice that $(x_n)$ is $\al^{1/2}$-equivalent to $(u_n)$ with respect to $\|\cdot\|_p$
and $|\cdot|_{p,r}$ norms. Picking $\delta>0$ sufficiently small we ensure that the
sequence $(x_n)$ is $\al$-equivalent to $(v_{l_n})$ both in $\|\cdot\|_p$ and
$|\cdot|_{p,r}$ norms. Let $T:[x_n]\rightarrow [v_{l_n}]$ be the isomorphism defined by
$Tx_n=v_{l_n}$, $n\in\N$.

Put $z_n=Ty_n$, $n\in\N$. By Lemma \ref{stab}, using $\al$-equivalence of $(y_n)$ and
$(z_n)$, we obtain for any $a_1,\dots,a_N\in\R$ the following inequalities:
$$\al^{-4p-2}\|a_1z_1+\dots+a_Nz_N\|_p\le \sqrt[p]{M}|
a_1z_1+\dots+a_Nz_N |_{p,r}\le 6(p+q)\al^6 \|a_1z_1+\dots+a_Nz_N\|_p$$

Therefore by Remark \ref{al} the subspace $Y=[z_n]$ of $X$ satisfies the desired
condition. \hfill $\square$

\end{document}